\documentstyle[12pt,fullpage]{article}

\newcommand{\Z}{{\rm Z\kern-.35em Z}}
\newcommand{\bP}{{\rm |\kern-.15em P}}
\newcommand{\Q}{\kern.3em\rule{.07em}{.65em}\kern-.3em{\rm Q}}
\newcommand{\R}{{\rm I\kern-.15em R}}
\newcommand{\D}{{\rm |\kern-.15em D}}
\newcommand{\h}{{\rm |\kern-.15em H}}
\newcommand{\C}{\kern.3em\rule{.07em}{.65em}\kern-.3em{\rm C}}
\newcommand{\T}{{\rm T\kern-.35em T}}

\begin{document}
\title{A geometric inequality for circle packings}
\author{Pablo A. Parrilo
\thanks{Institut f\"ur Automatik,
        Physikstrasse 3 / ETL,
        ETH-Zentrum,
        CH-8092 Z\"urich, Switzerland. Email: \texttt{parrilo@aut.ee.ethz.ch}}
\and    Ronen Peretz\thanks{
        Department of Mathematics,
        Ben Gurion University of the Negev,
        Beer-Sheva, 84105,
        Israel. Email: \texttt{ronenp@math.bgu.ac.il} }}

\maketitle
 
\newtheorem{definition}{Definition}
\newtheorem{remark}{Remark}
\newtheorem{example}{Example}
\newtheorem{lemma}{Lemma}
\newtheorem{proposition}{Proposition}
\newtheorem{corollary}{Corollary}
\newtheorem{theorem}{Theorem}

\abstract{ A geometric inequality among three triangles, originating
  in circle packing problems, is introduced.  In order to prove it, we
  reduce the original formulation to the nonnegativity of a particular
  polynomial in four real indeterminates. Techniques based on sum of
  squares decompositions, semidefinite programming, and symmetry
  reduction are then applied to provide an easily verifiable
  nonnegativity certificate.  }

\section{Introduction}

In this paper we prove the following geometric inequality: suppose
that we have three triangles. One with sides of lengths $X,Y$ and $Z$,
a second with sides of lengths $U,V$ and $W$ and a third triangle with
sides of lengths $(X+U), (Y+V)$ and $(Z+W)$. Let us denote by $\alpha
$ the angle in the first triangle between the sides of lengths $X$ and
$Y$. Let $\beta $ be the corresponding angle in the second triangle
between $U$ and $V$ and let $\gamma $ be the corresponding angle
between $(X+U)$ and $(Y+V)$ in the third triangle. Then
$$
\alpha\cdot (X+Y-Z)+\beta\cdot (U+V-W) \le \gamma\cdot ((X+U)+(Y+V)-(Z+W)).
$$
It turns out that proving this inequality is not at all simple. The
need for this inequality originates in \cite{p}. This last paper
describes a new approach to circle packings \cite{b,mr,rs}. The main
features of this approach are the theory of Perron-Frobenius for
non-negative matrices, \cite{bp,ho,s} and fixed-point theory,
\cite{or,be}. In particular our approach uses a converse of the
contraction principle as appears in \cite{be}. A central object that
is introduced in \cite{p} is {\bf the $\overline{a}$-mapping,
$f_{\overline{a}}\,:\,\R^{+|V|} \rightarrow\R^{+|V|}$, of a graph
embedding}. This is a variant of Thurston's relaxation mapping. A key
property of the $\overline{a}$-mapping, $f_{\overline{a}}$, is its
{\bf super-additivity},
$$
\forall\,\overline{r},\overline{s}\in\R^{+|V|},\,\,\,f_{\overline{a}}(\overline{r})+
f_{\overline{a}}(\overline{s})\le f_{\overline{a}}(\overline{r}+\overline{s}).
$$
It turns out that this property is implied by the above geometric
inequality.  A very interesting feature of our proof of this
inequality is the use of semidefinite programming based algorithms for
producing representations of non-negative polynomials as sums of
squares \cite{ChoiLamReznick,Parrilo}.

We make a reduction of the original inequality to the non-negativity
of a certain real polynomial of degree 20 in 4 real indeterminates.
Although instances of this size are near the limits of what can be
achieved using generic methods, our particular polynomial enjoys
certain convenient sparsity and symmetry properties.  These allow the
use of new algorithms, introduced in \cite{GatermannParrilo},
customized for polynomials with an invariant structure. Even though
the computational procedures in its current form use floating point
arithmetic to arrive at the result, the final solution can be easily
verified in a completely independent fashion. The methods work quite
nicely in the problem at hand, producing a concise representation of
the polynomial as a sum of five squares, thereby concluding the proof.

\section{The super-additivity of $f_{\overline{a}}$}

As proved in \cite{p}, the super-additivity of $f_{\overline{a}}$ follows from 
\begin{theorem}
If $a,b,c,d,R,S > 0$, then
$$
R\sin^{-1}\left\{\sqrt{\frac{ab}{(R+a)(R+b)}}\right\}+
S\sin^{-1}\left\{\sqrt{\frac{cd}{(S+c)(S+d)}}\right\} \le
$$
$$
\le (R+S)\sin^{-1}\left\{\sqrt{\frac{(a+c)(b+d)}{(R+S+a+c)(R+S+b+d)}}\right\}.
$$
\end{theorem}
The theorem has two simple geometric interpretations: \\
\\
{\bf (1)} Three circles of radii $R,a$ and $b$ that are mutually
tangent to one another from the outside form an Euclidean triangle.
The vertices of the triangle are the centers of the circles. The sides
of the triangle have the following lengths: $R+a$, $a+b$ and $R+b$.
Similarly three circles of radii $S,c$ and $d$ form a triangle of
sides $S+c$, $c+d$ and $S+d$. Finally, a third such triangle is formed
by three circles of radii $R+S+a+c$, $a+c+b+d$ and $R+S+b+d$. We note
that the sides of the third triangle have lengths which are the sums
of the corresponding sides of the first two triangles.  On the other
hand the three sets of triples of circles form also three circular
triangles.  The vertices of these triangles are the tangency points of
pairs of circles in each triple.  The lemma implies that the circular
sides of the third (largest) triangle are greater than
or equal to the sums of the corresponding circular sides of the first two circular triangles. \\
\\
{\bf (2)} Let us consider three Euclidean triangles. One with sides of
lengths $X,Y$ and $Z$, and an angle $\alpha $ between $X$ and $Y$. A
second triangle with sides of lengths $U,V$ and $W$, and an angle
$\beta $ between $U$ and $V$. A third triangle with sides of lengths
$(X+U),(Y+V)$ and $(Z+W)$, and an angle $\gamma $ between $(X+U)$ and
$(Y+V)$.  Then
$$
\alpha\cdot (X+Y-Z) + \beta\cdot (U+V-W) \le \gamma\cdot ((X+U)+(Y+V)-(Z+W)).
$$

\section{Reducing to algebraic inequalities}

We now make a reduction of the inequality of Theorem 1. There
are two ideas involved in it. The first idea is summarized in the following,

\begin{lemma}
Suppose that there exists a twice differentiable, surjective and strictly increasing
function $f\,:\,I\rightarrow [0,1]$ which satisfies the following two conditions: \\
\\
(1)
$$
f^{''}(1-f^{2})+f\cdot (f^{'})^{2}\le 0\,\,\,\,on\,\,I.
$$
(2)
$$
\left(\frac{R}{R+S}\right)f^{-1}\left(\sqrt{\frac{ab}{(R+a)(R+b)}}\right)+
\left(\frac{S}{R+S}\right)f^{-1}\left(\sqrt{\frac{cd}{(S+c)(S+d)}}\right)\le
$$
$$
\le f^{-1}\left(\sqrt{\frac{(a+c)(b+d)}{(R+S+a+c)(R+S+b+d)}}\right),
$$
for all $a,b,c,d,R,S > 0$. Then, the inequality of Theorem 1 holds true.
\end{lemma}
{\bf Proof.} \\
Consider the function $y=\sin^{-1} f(x)$ for $x\in I$. Then,
$$
\frac{dy}{dx}=\frac{f^{'}}{\sqrt{1-f^{2}}},
$$
$$
\frac{d^{2}y}{dx^{2}}=\frac{f^{''}(1-f^{2})+f\cdot (f^{'})^{2}}{(1-f^{2})^{3/2}}.
$$
By condition (1) we get $d^{2}y/dx^{2}\le 0$ on $I$ and hence $y$ is concave in $I$.
So for any $x,z\in I$ and for any $0\le t\le 1$, we have
\begin{equation}
t\sin^{-1} f(x)+(1-t)\sin^{-1} f(z)\le\sin^{-1} f(tx+(1-t)z).
\end{equation}
We make the following choice,
$$
x=f^{-1}\left(\sqrt{\frac{ab}{(R+a)(R+b)}}\right),\,\,
z=f^{-1}\left(\sqrt{\frac{cd}{(S+c)(S+d)}}\right),\,\,t=\left(\frac{R}{R+S}\right).
$$
Then, by inequality (1) we get,
\begin{eqnarray}
\left(\frac{R}{R+S}\right) & \sin^{-1} & \left(\sqrt{\frac{ab}{(R+a)(R+b)}}\right)+ \\
+\left(\frac{S}{R+S}\right) & \sin^{-1} & \left(\sqrt{\frac{cd}{(S+c)(S+d)}}\right)\le
\sin^{-1} f(tx+(1-t)z) \nonumber .
\end{eqnarray}
By condition (2) we have,
$$
tx+(1-t)z\le f^{-1}\left(\sqrt{\frac{(a+c)(b+d)}{(R+S+a+c)(R+S+b+d)}}\right)
$$
and since $f$ is increasing and also $\sin^{-1}$ is increasing, we get,
\begin{equation}
\sin^{-1} f(tx+(1-t)y)\le\sin^{-1}\left(\sqrt{\frac{(a+c)(b+d)}{(R+S+a+c)(R+S+b+d)}}\right).
\end{equation}
Theorem 1 follows by inequalities (2) and (3). $\Diamond $ \\
\\
{\bf special cases.} \\
{\bf (I)} $f(x)=\sin x,\,\,I=[0,\pi/2]$. Then in this case we have,
$$
f^{''}(1-f^{2})+f\cdot (f^{'})^{2}=-\sin x\cos^{2} x + \sin x\cos^{2} x\equiv 0
$$
and condition (1) of the theorem is satisfied. Condition (2) is the inequality of Theorem 1
and so the theorem is correct trivially in this case. \\
\\
{\bf (II)} $f(x)=1-1/x,\,\,I=[1,\infty]$. Then in this case we have,
$$
f^{''}(1-f^{2})+f\cdot (f^{'})^{2}=-\frac{2}{x^{3}}\left(1-\left(1-\frac{1}{x}\right)^{2}\right)+
\left(1-\frac{1}{x}\right)\cdot\frac{1}{x^{4}}=
$$
$$
=\frac{1}{x^{5}}-\frac{3}{x^{4}} < 0,
$$
for $x\ge 1$. So condition (1) is satisfied. Condition (2) and the conclusion of
the theorem prove the following,

\begin{lemma}
If for every $a,b,c,d,R,S > 0$ the following inequality is true,
\begin{eqnarray*}
& & \left(\frac{R}{R+S}\right)\left(\frac{1}{1-\sqrt{(ab)/[(R+a)(R+b)]}}\right)+ \\
& + & \left(\frac{S}{R+S}\right)\left(\frac{1}{1-\sqrt{(cd)/[(S+c)(S+d)]}}\right)\le \\
& & \le\frac{1}{1-\sqrt{[(a+c)(b+d)]/[(R+S+a+c)(R+S+b+d)]}},
\end{eqnarray*}
then, the inequality of Theorem 1 is true.
\end{lemma}
The second idea in this approach (after that of Lemma 1) is an elementary trick to
get rid of the square root functions in Lemma 2. Let us denote,
$$
\alpha=\sqrt{\frac{a}{R+a}},\,\,\beta=\sqrt{\frac{b}{R+b}},\,\,
\gamma=\sqrt{\frac{c}{S+c}},\,\,\delta=\sqrt{\frac{d}{S+d}}.
$$
Then $0\le\alpha,\beta,\gamma,\delta\le 1$. Also $\alpha,\,\beta$ are independent
except for $\alpha=1$ iff $\beta=1$. That happens only if $R=0$. $\gamma,\,\delta$
are independent except for $\gamma=1$ iff $\delta=1$. That happens only if $S=0$.
For the inverse transformations we have,
$$
a=\left(\frac{\alpha^{2}}{1-\alpha^{2}}\right)R,\,b=\left(\frac{\beta^{2}}{1-\beta^{2}}\right)R,\,
c=\left(\frac{\gamma^{2}}{1-\gamma^{2}}\right)S,\,d=\left(\frac{\delta^{2}}{1-\delta^{2}}\right)S.
$$
With these notations, the left hand side of the inequality in Lemma 2 is,
$$
\left(\frac{R}{R+S}\right)\left(\frac{1}{1-\alpha\beta}\right)+
\left(\frac{S}{R+S}\right)\left(\frac{1}{1-\gamma\delta}\right)=
\frac{R(1-\gamma\delta)+S(1-\alpha\beta)}{(R+S)(1-\alpha\beta)(1-\gamma\delta)}.
$$
As for the right hand side, we have,
$$
I_{1}=\sqrt{\frac{a+c}{R+S+a+c}}=\sqrt{\frac{\alpha^{2}(1-\gamma^{2})R+\gamma^{2}(1-\alpha^{2})S}
{(1-\gamma^{2})R+(1-\alpha^{2})S}},
$$
$$
I_{2}=\sqrt{\frac{b+d}{R+S+b+d}}=\sqrt{\frac{\beta^{2}(1-\delta^{2})R+\delta^{2}(1-\beta^{2})S}
{(1-\delta^{2})R+(1-\beta^{2})S}}.
$$
Plugging these into the inequality of Lemma 2 we get,
$$
\frac{R(1-\gamma\delta)+S(1-\alpha\beta)}{(R+S)(1-\alpha\beta)(1-\gamma\delta)}\le
\frac{1}{1-I_{1}I_{2}}.
$$
Hence,
$$
I_{1}I_{2}\ge\frac{R\alpha\beta(1-\gamma\delta)+S\gamma\delta(1-\alpha\beta)}
{R(1-\gamma\delta)+S(1-\alpha\beta)}.
$$
On squaring both sides we conclude that in order to prove Theorem 1,
it suffices to prove the following,

\begin{lemma}
If $R,S > 0$ and $0 < \alpha,\beta,\gamma,\delta < 1$, then,
$$
\left(\frac{\alpha^{2}(1-\gamma^{2})R+\gamma^{2}(1-\alpha^{2})S}
{(1-\gamma^{2})R+(1-\alpha^{2})S}\right)
\left(\frac{\beta^{2}(1-\delta^{2})R+\delta^{2}(1-\beta^{2})S}
{(1-\delta^{2})R+(1-\beta^{2})S}\right)\ge
$$
$$
\ge\left(\frac{R\alpha\beta(1-\gamma\delta)+S\gamma\delta(1-\alpha\beta)}
{R(1-\gamma\delta)+S(1-\alpha\beta)}\right)^{2}.
$$
\end{lemma}
{\bf Proof:}
Let us define,
$$
E=\left(\frac{\alpha^{2}(1-\gamma^{2})R+\gamma^{2}(1-\alpha^{2})S}
{(1-\gamma^{2})R+(1-\alpha^{2})S}\right)
\left(\frac{\beta^{2}(1-\delta^{2})R+\gamma^{2}(1-\alpha^{2})S}
{(1-\delta^{2})R+(1-\alpha^{2})S}\right)-
$$
$$
-\left(\frac{R\alpha\beta(1-\gamma\delta)+S\gamma\delta(1-\alpha\beta)}
{R(1-\gamma\delta)+S(1-\alpha\beta)}\right)^{2}.
$$
Then,
$$
E=\frac{RS(R+S)[(1-\gamma\delta)L\cdot R+(1-\alpha\beta)M\cdot S]}
{[(1-\gamma^{2})R+(1-\alpha^{2})S][(1-\delta^{2})R+(1-\alpha^{2})S][(1-\gamma\delta)R
+(1-\alpha\beta)S]^{2}},
$$
where we have,
$$
L = \alpha^{2}\beta^{2}(\alpha-\beta)^{2}+(\alpha-\beta)^{2}\gamma^{3}\delta^{3}+
$$
$$
+ \{(\alpha\delta)^{2}(1-\alpha\beta)(1+\alpha\beta-2\beta^{2})-
(\alpha\delta)(\beta\gamma)(2-4\alpha\beta+\beta\alpha^{3}+\alpha\beta^{3})+
$$
$$
+(\beta\gamma)^{2}(1-\alpha\beta)(1+\alpha\beta-2\alpha^{2})\} + 
$$
$$
+ \{\gamma^{2}\beta(1-\alpha\beta)(2\alpha-\beta-\alpha\beta^{2})-
\gamma\delta(\alpha^{2}+\beta^{2}+2\alpha^{3}\beta^{3}-4\alpha^{2}\beta^{2})+
$$
$$
+ \delta^{2}\alpha(1-\alpha\beta)(2\beta-\alpha-\alpha^{2}\beta)\}\gamma\delta,
$$
or as a polynomial in $\gamma $ and $\delta $,
$$
L=\alpha^{2}\beta^{2}(\alpha-\beta)^{2}+\beta^{2}(1-\alpha\beta)(1+\alpha\beta-2\alpha^{2})\gamma^{2}+
$$
$$
+\alpha^{2}(1-\alpha\beta)(1+\alpha\beta-2\beta^{2})\delta^{2}-
\alpha\beta(2+\alpha\beta^{3}-4\alpha\beta+\beta\alpha^{3})\gamma\delta+
$$
$$
+\beta(1-\alpha\beta)(2\alpha-\beta-\alpha\beta^{2})\gamma^{3}\delta-
(\alpha^{2}+\beta^{2}+2\alpha^{3}\beta^{3}-4\alpha^{2}\beta^{2})\gamma^{2}\delta^{2}+
$$
$$
+\alpha(1-\alpha\beta)(2\beta-\alpha-\alpha^{2}\beta)\gamma\delta^{3}+
(\alpha-\beta)^{2}\gamma^{3}\delta^{3},
$$
and where
$M=M(\alpha,\beta,\gamma,\delta)=L(\gamma,\delta,\alpha,\beta)$. Thus
it suffices to prove that for any $0 < \alpha,\beta,\gamma,\delta < 1$
we have $L(\alpha,\beta,\gamma,\delta)\ge 0$. For this will also imply
that $M(\alpha,\beta,\gamma,\delta)\ge 0$ for any such a choice. This,
in turn, will show that $E\ge 0$ for every choice of $R,S > 0$ and $0
< \alpha,\beta,\gamma,\delta < 1$ and hence will prove Lemma 3. To
check the non-negativity of $L(\alpha,\beta,\gamma,\delta)$ we make
the substitutions
\begin{equation}
(\alpha,\beta,\gamma,\delta)=(\frac{x^{2}}{1+x^{2}},\frac{y^{2}}{1+y^{2}},
\frac{z^{2}}{1+z^{2}},\frac{w^{2}}{1+w^{2}})
\label{eq:cov}
\end{equation}
and clear the denominators. This will give us a polynomial in
$\R[x,y,z,w]$. In fact, this polynomial is
$$
P(x,y,z,w)=
$$
\begin{equation}
=L\left(
\frac{x^{2}}{1+x^{2}},\frac{y^{2}}{1+y^{2}},
\frac{z^{2}}{1+z^{2}},\frac{w^{2}}{1+w^{2}}
\right)
(1+x^{2})^{4}(1+y^{2})^{4}(1+z^{2})^{3}
(1+w^{2})^{3}.
\label{eq:ourpoly}
\end{equation}
As a consequence, it suffices to check that $P(x,y,z,w)$ is
non-negative for all real values of its indeterminates. We conclude
the proof of Lemma~3 in the following section, after a brief detour
explaining the sum of squares based methods we have used.

\section{Sums of squares}

An obvious sufficient condition for non-negativity is to represent
$P(x,y,z,w)$ as a sum of squares of real polynomials. The connections
between sums of squares and non-negativity have been extensively
studied since the end of the 19th century, when Hilbert showed that in
the general case the two conditions are not equivalent. We refer the
reader to the wonderful survey \cite{Reznick} by Reznick on the
available results and history of Hilbert's 17th problem. In the work
of Choi, Lam, and Reznick \cite{ChoiLamReznick} the algebraic
structure of sums of squares decompositions is fully analyzed, and the
important ``Gram matrix'' method is introduced.  On the computational
side, convex optimization approaches to this problem originate in the
early work of Shor \cite{Shor}. Recently, efficient techniques using
semidefinite programming and exploiting problem structure have been
developed in \cite{Parrilo,GatermannParrilo}. A brief description of
the methods follows, referring the reader to the cited works, and the
references therein, for the full algorithmic details.

We explain next the general idea of the Gram matrix method. Given a
multivariate polynomial $F(\mathbf{x})$ for which we want to decide
whether a sum of squares decomposition exists, we attempt to express
it as a quadratic form in a new set of variables $\mathbf{u}$. A
judicious choice of these new variables will depend on both the
sparsity structure and symmetry properties of $F$
\cite{ChoiLamReznick,GatermannParrilo}.  For instance, for the
simplest case of a generic dense polynomial of total degree $2 d$, the
variables $\mathbf{u}$ will be all the monomials (in the variables
$\mathbf{x}$) of degree less than or equal to $d$.  Consequently, we
try to represent $F(\mathbf{x})$ as:
\begin{equation}
F(\mathbf{x}) = \mathbf{u}^T Q \mathbf{u} 
\label{sosrep}
\end{equation}
where $Q$ is a constant matrix.  Since in general the variables
$\mathbf{u}$ will not be algebraically independent, the matrix $Q$ in
the representation (\ref{sosrep}) \emph{is not unique}. In fact, there
is an affine subspace of matrices $Q$ that satisfy the equality, as
can be easily seen by expanding the right-hand side and equating term
by term.  If in the representation above the matrix $Q$ can be chosen
to be \emph{positive semidefinite}, then a factorization of the matrix
$Q$ directly provides a sum of squares decomposition of
$F(\mathbf{x})$.  Conversely, if $F$ is a sum of squares, then such a
$Q$ can always be constructed by expanding the terms in monomials.
Therefore, the problem of checking if a polynomial can be decomposed
as a sum of squares is \emph{equivalent} to verifying whether a
certain affine matrix subspace intersects the cone of positive
definite matrices.  This latter class of convex optimization problems
is known as \emph{semidefinite programs} (SDP) \cite{HandSDP}, and can
be efficiently solved using a variety of numerical algorithms, mainly
based on interior point methods.

\begin{example}
  Consider the quartic form in two variables described below, and
  define ${\mathbf{u}} = [ x^2, y^2, x y]^T$.
\begin{eqnarray*}
F(x,y) &=& 2 x^4 + 2 x^3 y  - x^2 y^2 + 5 y^4 \\
&=&
\left[\begin{array}{c}
x^2 \\  y^2 \\ x y
\end{array}\right]^T
\left[\begin{array}{ccc}
q_{11} & q_{12} & q_{13} \\
q_{12} & q_{22} & q_{23} \\
q_{13} & q_{23} & q_{33}
\end{array}\right]
\left[\begin{array}{c}
x^2 \\  y^2 \\ x y
\end{array}\right]\\
&=&
q_{11} x^4 + q_{22} y^4 + (q_{33} + 2 q_{12}) x^2 y^2
+ 2 q_{13} x^3 y + 2 q_{23} x y^3
\end{eqnarray*}
Therefore, in order to have an identity, the following linear
equalities should hold:
\begin{equation}
q_{11} = 2, \quad
q_{22} = 5, \quad
q_{33} + 2 q_{12} = -1, \quad
2 q_{13} = 2, \quad
2 q_{23} = 0.
\end{equation}
A positive semidefinite $Q$ that satisfies the linear equalities can
then be found using semidefinite programming. A particular solution is
given by:
\[
Q =
\left[\begin{array}{rrr}
2  & -3 & 1 \\ -3 & 5 & 0 \\ 1 & 0 & 5
\end{array}\right]
= L^T L, \qquad
L =
\frac{1}{\sqrt{2}}\left[\begin{array}{rrr}
2 & -3 & 1 \\
0 & 1 & 3
\end{array}\right],
\]
and therefore we have the sum of squares decomposition:
\[
F(x,y) = \frac{1}{2} (2 x^2 - 3 y^2 + x y)^2 +
\frac{1}{2}(y^2 + 3 x y)^2.
\]
\label{ex:sosexample}
\end{example}

Back to our concrete problem, the polynomial $P$ in (\ref{eq:ourpoly})
has four variables, 123 nonzero monomial terms, and total degree 20.
Notice that a polynomial of that degree and number of variables
generically has ${20 + 4 \choose 4} = 10626$ monomials, so $P$ is
quite sparse. In particular, it has a bipartite structure, with degree
12 in $x,y$, and degree 8 in $z,w$. Furthermore, $P$ has very
appealing symmetry properties: some inherited from $L$, and some as a
result of the substitution (\ref{eq:cov}). Concretely, it is easy to
see that $P$ is invariant under the transformations:
\begin{eqnarray}
(x,y,z,w) & \rightarrow & (y,x,w,z)  \\
          & \rightarrow & (\pm x,\pm y,\pm z,\pm w)
\end{eqnarray}
The first property is a clear consequence of the symmetry of our
original geometric inequality with respect to interchange of the two
triangles. The second one is a side effect of our choice for modeling
the nonnegativity constraints. The transformations given above
generate a symmetry group $G$ with 32 elements and 14 irreducible
representations: eight one-dimensional and six two-dimensional ($8
\cdot 1^2 + 6 \cdot 2^2 = 32$).  As explained extensively by Gatermann
and Parrilo in \cite{GatermannParrilo}, these symmetries can be
exploited very successfully in reducing the computational
requirements.

To do this, the approach in \cite{GatermannParrilo} relies on a
crucial property of convex optimization problems invariant under a
group action, namely the fact that the optimal solution can always be
restricted to the fixed-point subspace. Using Schur's lemma of
representation theory, it is shown there that by using an appropriate
symmetry-adapted coordinate transformation, the original semidefinite
program can be decomposed into a collection of smaller coupled
problems, of cardinality equal to the number of irreducible
representations of the group. This reduces both the size and the
number of variables in the problem, and as a consequence notably
enhances both the accuracy and conditioning of the solution.

Attempting to directly establish the nonnegativity of $P$ without
taking into account both the sparsity and symmetries can be a
difficult (or even impossible) task for current SDP solvers, both in
terms of memory requirements and accuracy. A naive approach, using
only degree information but no structure whatsoever, would require
solving a semidefinite program of dimension $1001\times 1001$ and
10626 constraints.  By exploiting only the sparsity of $P$, but not
its symmetry, the problem is reduced to dimension $137\times 137$ and
1328 constraints.  Adding the simplifications resulting from the
symmetries, the problem is further simplified to a much more
manageable one with 14 coupled SDP (one for each irreducible
representation), of dimensions ranging between $2\times 2$ and $11
\times 11$ (see Table~\ref{tab:sdpdims}).  For instance, for the
trivial irreducible representation (\# 1 in the table), the
corresponding new variables $\mathbf{u}$ are invariant under the group
action, and given by:
\[
\begin{array}{c}
y^2z^2+x^2w^2 \\
x^2z^2w^2+y^2z^2w^2, \quad y^4z^2+x^4w^2, \quad x^2y^2z^2+x^2y^2w^2,
\quad x^4y^2+x^2y^4\\
x^2y^2z^2w^2, \quad x^4z^2w^2+y^4z^2w^2, \quad x^2y^4z^2+x^4y^2w^2,
\quad x^4y^2z^2+x^2y^4w^2.
\end{array}
\]
Notice in Table~\ref{tab:sdpdims} that the combined total, taking into
account multiplicities, is equal to 137, the dimension of the sparse
version of the problem.

\begin{table}
\begin{center}
\begin{tabular}{c|rrr rrr rrr rrr rr}
Irr. Rep. \# & 1 & 2 & 3 & 4 & 5 & 6 & 7 & 8 &9 &10 &11 &12 &13 &14 \\ \hline
Multiplicity&  1 & 1 & 1 & 1 & 1 & 1 & 1 & 1 &2 &2  &2  &2  &2  &2 \\
Dim. SDP    &  9 & 6 & 6 & 4 & 8 & 5 & 3 & 2 & 11 & 7 & 8 & 7 & 8 & 6
\end{tabular}
\end{center}
\caption{Irreducible representations of $G$ and the corresponding SDP dimensions.}
\label{tab:sdpdims}
\end{table}

The resulting system of matrix inequalities can be solved with
standard SDP solvers, such as SeDuMi \cite{sedumi}. The output
provides a decomposition of $P$ as a sum of squares of polynomials,
with coefficients given by floating point numbers. In this particular
case, the computed values immediately suggest the existence of a
solution, presented below, with polynomials having integer
coefficients. The solution can be verified in a completely
independent fashion, providing a mathematically correct certificate of
the nonnegativity of the polynomial $P$.

\vspace{.5 cm}
{\bf Proof of Lemma 3: (continued) }
Consider the following three polynomials:
$$
A(x,y,z,w)=-y^{2}z^{2}-y^{4}z^{2}+x^{2}w^{2}+2x^{2}y^{2}w^{2}-2x^{2}y^{2}z^{2}-x^{2}y^{4}-
$$
$$
-2x^{2}y^{4}z^{2}+x^{4}w^{2}+x^{4}y^{2}+2x^{4}y^{2}w^{2},
$$
$$
B(x,y,z,w)=(1+x^{2}+y^{2})(-x^{2}w^{2}-x^{2}z^{2}w^{2}-x^{2}y^{2}w^{2}+x^{2}y^{2}z^{2}+
$$
$$
+y^{2}z^{2}+y^{2}z^{2}w^{2})
$$
and
$$
C(x,y,z,w)=(x-y)(x+y)(-x^{2}z^{2}w^{2}+x^{2}y^{2}+x^{2}y^{2}w^{2}+x^{2}y^{2}z^{2}-
$$
$$
-z^{2}w^{2}-y^{2}z^{2}w^{2}).
$$
Then we have the following identity
$$
P(x,y,z,w)=A(x,y,z,w)^{2}(z^{2}+w^{2}+2z^{2}w^{2})+B(x,y,z,w)^{2}+C(x,y,z,w)^{2}.
$$
Thus $P(x,y,z,w)$ is a sum of five squares of real polynomials and the
proof of Lemma 3 is completed. $\Diamond $ \\ \\

Rewriting the obtained sum of squares decomposition in terms of the
original variables, the following representation of $L$ can be
obtained:
\begin{eqnarray*}
L &=&  L_1 + L_2 + L_3 \\
L_1 &=&
(\gamma+\delta) (-\alpha^2 \beta+\alpha \beta^2-\alpha \delta+\beta \gamma-\beta \gamma \delta+\alpha \delta \gamma-\alpha \beta^2 \gamma+\alpha^2 \beta \delta)^2 \\
L_2 &=& (1-\gamma) (1-\delta) 
(\alpha \beta-1)^2 (\alpha \delta-\beta \gamma)^2  \\
L_3 &=& (1-\gamma) (1-\delta) 
(\alpha-\beta)^2 (\alpha \beta-\gamma \delta)^2. 
\end{eqnarray*}
From this, stronger conclusions on the sign of $L$ can be immediately
derived: not only it is nonnegative on the open unit hypercube
$(0,1)^4$ as needed for Lemma~3, but the same property holds on the
much larger region $\R \times \R \times \{\gamma+\delta \geq 0,
(1-\gamma)(1-\delta) \geq 0\}$.

\vspace{1cm}
\noindent \textbf{Acknowledgments:} We would like to warmly
acknowledge the help of Bruce Reznick, who made possible this
collaboration by introducing the authors to each other's work.

\end{document}